\documentclass{amsart}
\newtheorem{theorem}{Theorem}
\newtheorem{lemma}[theorem]{Lemma}
\newcommand{\cc}{\mathbb{C}}

\newcommand{\dbar}{\overline{\partial}}
\newcommand{\ybar}{\overline{y}}
\newcommand{\zbar}{\overline{z}}

\title[H\"older estimates for the $\overline\partial$-equation]
{H\"older estimates for the $\overline\partial$-equation on 
surfaces with simple singularities }

\author{F. Acosta}
\author{E. S. Zeron}
\address{Depto. Matem\'aticas, CINVESTAV, Apartado 
Postal 14-740, M\'exico D.F., 07000, M\'exico.}
\email{facosta@math.cinvestav.mx}
\email{eszeron@math.cinvestav.mx}

\date{\today}
\thanks{Research supported by Cinvestav(Mexico) and Conacyt(Mexico)}

\begin{document}
\begin{abstract}
Let $\Sigma\subset\cc^3$ be a $2$-dimensional subvariety with an 
isolated simple (rational double point) singularity at the origin. 
The main objective of this paper is to solve the $\dbar$-equation 
on a neighbourhood of the origin in $\Sigma$, demanding a H\"older 
condition on the solution. 
\end{abstract}

\maketitle
\section{Introduction}

Let $\Sigma\subset\cc^3$ be a subvariety with an isolated singularity at 
the origin. Given a $\dbar$-closed $(0,1)$-differential form $\lambda$ 
defined on $\Sigma$ minus the origin, Gavosto and Forn{\ae}ss proposed a 
general technique for solving the differential equation $\dbar{g}=\lambda$ 
on a neighbourhood of the origin in $\Sigma$. The calculations were done 
in the sense of distributions, and they demanded an extra H\"older 
condition on the solution $g$, see \cite{Ga} and \cite{GF}. Their basic 
idea was to analyse $\Sigma$ as a branched covering over $\cc^2$, to 
solve the corresponding $\dbar$-equation on $\cc^2$, and to \textit{lift} 
the solution from $\cc^2$ into $\Sigma$ again. Gavosto and Forn{\ae}ss 
completed all the calculations in the particular case when 
$\Sigma\subset\cc^3$ is defined by the polynomial $x_1x_2=x_3^2$. That is, 
when $\Sigma$ is a surface with an isolated simple (rational double point) 
singularity of type $A_2$ at the origin, \cite[p.~60]{Di}.

Let $X_N$ and $Y_N$ be two subvarieties of $\cc^3$ defined by the 
respective polynomials $x_1x_2=x_3^N$ and $y_1^2y_3+y_2^2=y_3^{N+1}$, 
for any natural number $N\geq2$. Surface $X_N$ (respect. $Y_N$) has an 
isolated simple singularity of type $A_{N-1}$ (respect. $D_{N+2}$) at the 
origin, see \cite[p.~60]{Di}. The main objective of this paper is to give 
an alternative and simplified solution to the equation $\dbar{g}=\lambda$ 
on both surfaces $X_N$ and $Y_N$, with an extra H\"older condition on $g$. 
The central idea is to consider $\cc^2$ as a branched covering over $X_N$ 
and $Y_N$, instead of analysing $X_N$ as a branched covering over $\cc^2$. 
In the case of $X_N$, we use the natural branched $N$-covering 
$\pi_N:\cc^2\to{X_N}$ defined by $\pi_N(z_1,z_2)=(z_1^N,z_2^N,z_1z_2)$, 
in order to obtain the following theorem. We shall explain, at the end of 
the third section of this paper, why we use the covering $\pi_N$ instead 
of a standard \textit{blow up} mapping. 

\begin{theorem}\label{main}
Let $Ev(N)$ be the smallest even integer greater than or equal to $N$. 
Given an open ball $B_R\subset\cc^2$ of radius $R>0$ and centre in the 
origin, there exists a finite positive constant $C_1(R)$ such that: For 
every continuous $(0,1)$-differential form $\lambda$ defined on the 
compact set $\pi_N(\overline{B_R})\subset{X_N}$, and $\dbar$-closed 
on the interior $\pi_N(B_R)$, the equation $\dbar{h}=\lambda$ has a 
continuous solution $h$ on $\pi_N(B_R)$ which also satisfies the 
following H\"older estimate, with $\beta=1/Ev(N)$,
\begin{equation}\label{eq1} 
\|h\|_{\pi_N(B_R)}+\sup_{x,w\in\pi_N(B_R)}\frac{|h(x)-h(w)|}
{\|x-w\|^{\beta}}\,\leq{C}_1(R)\|\lambda\|_{\pi_N(B_R)}. 
\end{equation} 
\end{theorem}

In the last section of this paper, we extend Theorem~\ref{main} to 
solve the $\dbar$-equation on the subvariety $Y_N$ as well. The notation 
$\|h\|_K$ stands for the maximum of $|h|$ on the compact set $K$, and 
$\|x-w\|$ stands for the euclidean distance between $x$ and $w$. Since 
$\|x-w\|$ is less than or equal to the distance between $x$ and $w$ 
measured \textit{along} the surface $X_N$, we can assert that 
inequality~(\ref{eq1}) is indeed a H\"older estimate on $X_N$ itself. 
Finally, all differentials are defined in terms of distributions. For 
example, the fact that the continuous $(0,1)$-differential form $\lambda$ 
is $\dbar$-closed on $\pi_N(B_R)$ means that the integral~:
\begin{equation}\label{eq1a} 
\int_{\pi_N(B_R)}\lambda\wedge\dbar\sigma\,=0, 
\end{equation} 
for every smooth $(2,0)$-differential form $\sigma$ defined on 
$\pi_N(B_R)\setminus\{0\}$, such that both $\sigma$ and $\dbar\sigma$ 
extend continuously to the origin, and these extensions have both 
compact support inside $\pi_N(B_R)$.

The proof of Theorem~\ref{main} is presented in the following two 
sections. The next section is devoted to introducing all the basic ideas 
for the particular case when $N=2$. Moreover, in the third section of this 
paper, we shall use these ideas for solving the $\dbar$-equation on $X_N$, 
in the extended case $N\geq3$. Finally, in the last section of this paper, 
we extend Theorem~\ref{main} to solve the $\dbar$-equation on the 
subvariety $Y_N$ as well

\section{Proof of Theorem~\ref{main}, case $N=2$.}

Consider the natural branched covering 
$\pi_2(z_1,z_2)=(z_1^2,z_2^2,z_1z_2)$ defined from $\cc^2$ onto 
$X_2:=[x_1x_2=x_3^2]$. It is easy to see that $\pi_2$ is a branched 
$2$-covering, and that the origin is the only branch point of $\pi_2$, 
because the inverse image $\pi_2^{-1}(x)$ is a set of the form 
$\{\pm{z}\}$, for every $x\in{X_2}$. Besides, define the antipodal 
automorphism $\phi(z)=-z$ which allows us to jump between the different 
\textit{branches} of $\pi_2$. In particular, we have that 
$\phi^*\pi_2(z)=\pi_2(-z)=\pi_2(z)$.

We assert that the operators $\pi_2^*$ and $\dbar$ commute. It is easy 
to see that $\pi_2^*$ and $\dbar$ commute when $\dbar$ is a standard 
differential, for $\pi_2$ is holomorphic. However, calculations become 
more complicated when $\dbar$ is analysed in the sense of distributions. 
Let $B_R\subset\cc^2$ be an open ball of radius $R>0$. We prove the 
commutativity of $\pi_2^*$ and $\dbar$ for the particular case of a 
$\dbar$-closed $(0,1)$-differential form $\lambda$ defined on 
$\pi_2(B_R)$; the proof with a general differential form follows 
exactly the same procedure. We have that $\dbar\lambda=0$ in the 
sense of equation~(\ref{eq1a}), and we need to prove that 
$\dbar(\pi_2^*\lambda)$ is equal to $\pi_2^*(\dbar\lambda)=0$ 
in the sense of distributions, that is~:
\begin{equation}\label{eq1b}
\int_{B_R}\pi_2^*\lambda\wedge\dbar{v}\,=0,
\end{equation}
for every smooth $(2,0)$-differential form $v$ with compact support 
in $B_R$. The automorphism $\phi$ preserves the orientation of $B_R$, 
for it is analytic. Thus, after doing a simple change of variables, 
and recalling that $\phi^*\pi_2=\pi_2$, we have that the integral 
in equation~(\ref{eq1b}) is equal to 
$\int_{B_R}\pi_2^*\lambda\wedge\dbar\phi^*v$. Moreover, since 
$v+\phi^*v$ is constant in the fibres of $\pi_2$ (it is invariant 
under the \textit{pull back} $\phi^*$) there exists a second 
differential form $\sigma$ defined on $\pi_2(B_R)$ such that 
$v+\phi^*v$ is equal to $\pi_2^*\sigma$. Hence~:
$$\int_{B_R}\pi_2^*\lambda\wedge\dbar{v}\,
=\int_{B_R}\pi_2^*\lambda\wedge\dbar\;\frac{v+\phi^*v}{2}\,
=\int_{\pi_2(B_R)}\frac{\lambda\wedge\dbar\sigma}{2}\,=0.$$
The equality to zero follows from equation~(\ref{eq1a}), and so 
$\dbar(\pi_2^*\lambda)=0$ on $B_R$, as we wanted to prove. Suppose now 
that the differential equation $\dbar{g}=\pi_2^*\lambda$ has a solution 
$g$ on $B_R$. The sum $g+\phi^*g$ is also constant in the fibres of 
$\pi_2$ (it is invariant under the \textit{pull back} $\phi^*$), so 
there exists a continuous function $f$ on $B_R$ such that $\pi_2^*f$ 
is equal to $g+\phi^*g$. We assert that $\dbar{f}=2\lambda$ on 
$\pi_2(B_R)$. This result follows automatically because~:
$$\pi_2^*\dbar{f}=\dbar(g+\phi^*g)=\pi_2^*\lambda
+\phi^*\pi_2^*\lambda=\pi_2^*(2\lambda).$$

Previous equation demands that the operators $\phi^*$ and $\dbar$ commute 
as well in $B_R$, when $\dbar$ is seen as a distribution. This is an 
exercise based on the fact that integral $\int\phi^*\aleph=\int\aleph$, 
as we have indicated in the paragraph situated after 
equation~(\ref{eq1b}), and because $\phi$ preserves the orientation of 
$B_R$. Suppose now that $\lambda$ is also continuous on the compact set 
$\pi_2(\overline{B_R})$. Then, we can apply Theorems~2.1.5 and~2.2.2 of 
\cite{HL} in order to get the following H\"older estimate.

\begin{theorem}\label{Henkin}
Given an open ball $B_R\subset\cc^2$ of radius $R>0$ and centre in the 
origin, there exist two finite positive constants $C_2(R)$ and $C_3(R)$ 
such that: For every continuous $(0,1)$-differential form $\lambda$ 
defined on $\pi_2(\overline{B_R})\subset{X}$, and $\dbar$-closed on 
the interior $\pi_2(B_R)$, the equation $\dbar{g}=\pi_2^*\lambda$ has 
a continuous solution $g$ on $B_R$ which also satisfies the following 
H\"older estimates,
\begin{eqnarray}\label{eq2}
\|g\|_{B_R}+\sup_{z,\zeta\in{B_R}}\frac{|g(z)-g(\zeta)|}
{\|z-\zeta\|^{1/2}}&\leq&C_2(R)\|\pi_2^*\lambda\|_{B_R},\\
\label{eq3}\hbox{and}\quad\sup_{z,\zeta\in{B_{R/2}}}\frac{|g(z)
-g(\zeta)|}{\|z-\zeta\|}&\leq&C_3(R)\|\pi_2^*\lambda\|_{B_R}.  
\end{eqnarray} 
\end{theorem}

\begin{proof} Inequality (\ref{eq2}) holds because of Theorem~2.2.2 in 
\cite{HL}. Besides, recalling the proofs of Lemma~2.2.1 and Theorem~2.2.2, 
in \cite{HL}, we have that inequality (\ref{eq3}) holds whenever there 
exists a finite positive constant $C_4(R)$ such that~:
\begin{equation}\label{eq4}
\sup_{z,\zeta\in{B}_{R/2}}\frac{|E(z)-E(\zeta)|}{\|z-\zeta\|}\,\leq
C_4(R)\|\pi_2^*\lambda\|_{B_R},
\end{equation}
for every function $E(z)$ defined according to equation~(2.2.7) of 
\cite[p.~70]{HL}. Let $\Upsilon$ be the closed interval which joins 
$z$ and $\zeta$ inside the ball $B_{R/2}$. Then,
\begin{eqnarray}\label{eq4a}|E(z)-E(\zeta)|
&\leq&\int_0^1\left|\frac{d}{dt}E(t\zeta+(1-t)z)\right|\,dt\\
\nonumber&\leq&\|z-\zeta\|\,\sup_{y\in\Upsilon}\,\sum^2_{k=1}
\left|\frac{\partial{E}}{\partial{y}_k}\right|+
\left|\frac{\partial{E}}{\partial\ybar_k}\right|.
\end{eqnarray}

Finally, by equation~(2.2.9) in \cite{HL}, we know there exists 
a finite constant $C_4(R)$ such that all partial derivatives 
$\left|\frac{\partial{E}}{\partial{y}_k}\right|$ and 
$\left|\frac{\partial{E}}{\partial\ybar_k}\right|$ are less than or 
equal to \linebreak $\frac{C_4(R)}{5}\|\pi_2^*\lambda\|_{B_R}$, for 
every $y\in{B}_{R/2}$ and each index $k=1,2$. Notice that $D=B_R$ in 
equations~(2.2.7) and~(2.2.9), but $y$ lies inside the smaller ball 
$B_{R/2}$. Thus, equation~(\ref{eq4a}) automatically implies that 
inequalities~(\ref{eq4}) and~(\ref{eq3}) holds, as we wanted. 
\end{proof}

The problem is now reduced to estimating the distance $\|z-\zeta\|$ 
with respect to the projections $\|\pi_2(z)-\pi_2(\zeta)\|$.

\begin{lemma}\label{estimate}
Given two points $z$ and $\zeta$ in $\cc^2$ such that $\|z-\zeta\|$ is 
less than or equal to $\|z+\zeta\|$, the following inequality holds.
$$2\|\pi_2(z)-\pi_2(\zeta)\|\,\geq\,\|z-\zeta\|
\max\{\|z\|,\|\zeta\|,\|z-\zeta\|\}.$$
\end{lemma}

\begin{proof} We know that $2\|z\|$ and $2\|\zeta\|$ are both less than 
or equal to $\|z+\zeta\|+\|z-\zeta\|$. The given hypotheses indicates that 
$\|z-\zeta\|\leq\|z+\zeta\|$. Hence, the maximum of $\|z\|$, $\|\zeta\|$ 
and $\|z-\zeta\|$ is also less than or equal to $\|z+\zeta\|$. The wanted 
result will follows after proving that $\|z-\zeta\|\cdot\|z+\zeta\|$ is 
less than or equal to $2\|\pi_2(z)-\pi_2(\zeta)\|$. Setting 
$P_1=z_1-\zeta_1$, $P_2=z_2-\zeta_2$, $Q_1=z_1+\zeta_1$ and 
$Q_2=z_2+\zeta_2$, allows us to write the following series 
of inequalities:
\begin{eqnarray*}
&&\|z-\zeta\|^2\cdot\|z+\zeta\|^2\;=\\
&=&|P_1Q_1|^2+|P_1Q_2|^2+|P_2Q_1|^2+|P_2Q_2|^2\\
&\leq&4|P_1Q_1|^2+4|P_2Q_2|^2+|P_1Q_2|^2+|P_2Q_1|^2-2|P_1Q_1P_2Q_2|\\
&\leq&4|P_1Q_1|^2+4|P_2Q_2|^2+|P_1Q_2+P_2Q_1|^2\\
&=&4\|\pi_2(z)-\pi_2(\zeta)\|^2.
\end{eqnarray*}
\end{proof}

We are now in position to prove Theorem~\ref{main} 
for the simplest case $N=2$.

\begin{proof}\textbf{(Theorem~\ref{main}, case $N=2$)}. 
Suppose that $\lambda=\sum\lambda_kd\overline{x}_k$. Then, 
\begin{equation}\label{eq5}
\pi_2^*\lambda=[2\zbar_1\lambda_1(\pi_2)+\zbar_2\lambda_3(\pi_2)]
d\zbar_1+[2\zbar_2\lambda_2(\pi_2)+\zbar_1\lambda_3(\pi_2)]d\zbar_2.
\end{equation}
We obviously have that $|z_k|<R$ for every point $z\in{B_R}$. Hence,
\begin{equation}\label{eq6}
\|\pi_2^*\lambda\|_{B_R}\leq3R\,\|\lambda\|_{\pi_2(B_R)}.
\end{equation}

Let $g$ be a continuous solution to the equation $\dbar{g}=\pi_2^*\lambda$
on $B_R$, and suppose that $g$ satisfies the H\"older estimates given in
equations~(\ref{eq2}) and (\ref{eq3}) of Theorem~\ref{Henkin}. Recalling 
the analysis done in the paragraphs situated before Theorem~\ref{Henkin}, 
we know there exists a continuous function $h$ defined on $\pi_2(B_R)$ 
such that $h\circ\pi_2$ is equal to $\frac{g+\phi^*g}{2}$. In particular, 
$\dbar{h}=\lambda$ on $\pi_2(B_R)$, and
\begin{equation}\label{eq7}
\|h\|_{\pi_2(B_R)}=\frac{\|g+\phi^*g\|_{B_R}}{2}\,\leq\|g\|_{B_R}.
\end{equation}

Notice that $\beta=1/2$ when $N=2$. Given two points $x,w\in\pi_2(B_R)$, 
choose $z,\zeta\in{B_R}$ such that $x=\pi_2(z)$ and $w=\pi_2(\zeta)$. 
Since $\pi_2(\zeta)=\pi_2(-\zeta)$, we can even choose $\zeta\in{B_R}$ 
so that $\|z-\zeta\|$ is less than or equal to $\|z+\zeta\|$. If $z$ 
and $\zeta$ are both inside the ball $B_{R/2}$, we may apply 
equation~(\ref{eq3}) of Theorem~\ref{Henkin}, and the inequality 
$2\|x-w\|\geq\|z-\zeta\|^2$ given in Lemma~\ref{estimate}, in 
order to get~:
\begin{eqnarray}\label{eq8}
\frac{|h(x)-h(w)|}{2^{1/2}\|x-w\|^{1/2}}&\leq&
\frac{|g(z)-g(\zeta)|+|g(-z)-g(-\zeta)|}{2\,\|z-\zeta\|}\\
\nonumber&\leq&C_3(R)\|\pi_2^*\lambda\|_{B_R}.
\end{eqnarray}

On the other hand, suppose, without lost of generality, that $z$ 
is not inside the ball $B_{R/2}$; that is $\|z\|\geq\frac{R}{2}$. 
Lemma~\ref{estimate} implies then that $\|x-w\|$ is greater than 
or equal to $\frac{R}{4}\|z-\zeta\|$. Whence, equation~(\ref{eq2}) 
automatically implies the following,
\begin{equation}\label{eq9}
\frac{|h(x)-h(w)|}{\|x-w\|^{1/2}}\leq\frac{2}{\sqrt{R}}\,
C_2(R)\|\pi_2^*\lambda\|_{B_R}.
\end{equation}

Finally, considering Theorem~\ref{Henkin} and equations~(\ref{eq6}) to 
(\ref{eq9}), we can deduce the existence of a bounded positive constant 
$C_1(R)$ such that equation~(\ref{eq1}) holds.
\end{proof}

We close this section with some observations about Theorem~\ref{main}. 
Firstly, the procedure presented in this section yields a continuous 
solution $h$ to the equation $\dbar{h}=\lambda$. Moreover, we are directly 
using the estimates given in \cite{HL}, but we may use any integration 
kernel which produces estimates similar to those presented in 
equations~(\ref{eq2}) and (\ref{eq3}) of Theorem~\ref{Henkin}.

On the other hand, the extension of Theorem~\ref{main} to considering 
a general subvariety $\Sigma$, with an isolated singularity, does not 
seem to be trivial. Theorem~\ref{main} demands the existence of a 
branched finite covering $\pi:W\to\Sigma$, where $W$ is a \textit{nice} 
non-singular manifold and the inverse image of the singular point is a 
singleton. It does not seem to be trivial to produce such a branched 
finite covering.

\section{Proof of Theorem~\ref{main}, case $N\geq3$.}

We analyse in this section the general case of the variety 
$X_N\subset\cc^3$ defined by $x_1x_2=x_3^N$, for any natural number 
$N\geq3$. Surface $X_N$ has an isolated simple singularity of type 
$A_{N-1}$ at the origin, \cite[p.~60]{Di}. Define the automorphisms 
$\phi_k:\cc^2\to\cc^2$, for each natural number $k$,
\begin{equation}\label{eq11}
\phi_k(z_1,z_2)=(\rho_N^kz_1,\rho_N^{-k}z_2)
\quad\hbox{where}\quad\rho_N=e^{2\pi{i}/N}.
\end{equation}

Consider the natural branched covering 
$\pi_N(z_1,z_2)=(z_1^N,z_2^N,z_1z_2)$ defined from $\cc^2$ onto $X_N$. 
It is easy to see that $\pi_N$ is a branched $N$-covering, and that the 
origin is the only branch point of $\pi_N$, because the inverse image 
$\pi^{-1}_N(x)$ is a set of the form $\{\phi_k(z)\}_{1\leq{k}\leq{N}}$, 
for every $x\in{X_N}$. Thus, the automorphisms $\phi_k$ allow us to jump 
between the different \textit{branches} of $\pi_N$. In particular, we have 
that $\pi_N=\phi_k^*\pi_N$ for every $k$. Besides, the operators $\pi_N^*$ 
and $\dbar$ commute, the proof is based on the same ideas presented at the 
beginning of section two.

Given a $\dbar$-closed $(0,1)$-differential form $\lambda$ defined on 
$X_N$, we obviously have that $\dbar(\pi_N^*\lambda)=0$. Suppose the 
differential equation $\dbar{g}=\pi_N^*\lambda$ has a solution $g$ in 
$\cc^2$. The sum $\frac{1}{N}\sum_{k=1}^N\phi_k^*g$ is constant in the 
fibres of $\pi_N$ (it is invariant under every \textit{pull back} 
$\phi_j^*$), so there exists a continuous function $h$ on $X_N$ such 
that $\pi_N^*h$ is equal to $\frac{1}{N}\sum\phi_k^*g$. We assert that 
$\dbar{h}=\lambda$ on $X_N$. This result follows automatically because,
$$\textstyle{\pi_N^*\dbar{h}=\frac{1}{N}\sum\dbar\phi_k^*g
=\frac{1}{N}\sum\phi_k^*\pi_N^*\lambda=\pi_N^*\lambda}.$$

Let $B_R\subset\cc^2$ be an open ball of radius $R>0$. If $\lambda$ is 
continuous on the compact set $\pi_N(\overline{B_R})$, then, we can apply 
Theorem~\ref{Henkin} in order to get a solution $h$ which satisfies a 
H\"older estimate on the ball $B_R$. Obviously, like in the second 
section of this paper, the central part of the proof is an estimate 
of the distance $\|z-\zeta\|$ with respect to the projections 
$\|\pi_N(z)-\pi_N(\zeta)\|$. This estimate is done in the next 
lemma. Given two points $z$ and $\zeta$ in $\cc^2$, notation 
$\|z,\zeta\|_\infty$ stands for the maximum of $|z_1|$, $|z_2|$, $|\zeta_1|$ 
and $|\zeta_2|$. Moreover, $\|z\|_{\infty}:=\|z,0\|_{\infty}$ as well.

\begin{lemma}\label{general}
Let $Ev(N)$ be the smallest even integer greater than or equal 
to $N$. Given two points $z$ and $\zeta$ in $\cc^2$ such that 
$\|z-\phi_k(\zeta)\|_\infty$ is greater than or equal to 
$\|z-\zeta\|_\infty$ for every automorphism $\phi_k$ defined 
in~(\ref{eq11}), the following inequality holds for $\delta$ 
equal to both $N$ and $Ev(N)/2$.
\begin{equation}\label{eq12}
\|\pi_N(z)-\pi_N(\zeta)\|\geq\min\left\{\frac{\|z-\zeta\|^2_\infty}
{12},\frac{\|z,\zeta\|^{N-\delta}_\infty\|z-\zeta\|^\delta_\infty}
{(8/3)^{N-\delta}\,2^\delta}\right\}.
\end{equation}
\end{lemma}

\begin{proof} Set $z=(a,b)$, so that $\pi_N(z)=(a^N,b^N,ab)$. Moreover, 
given $\zeta=(s,t)$, we can suppose without loss of generality that 
$|a-s|\geq|b-t|$, and so $\|z-\zeta\|_\infty=|a-s|$. We shall prove 
inequality~(\ref{eq12}) by considering three cases.

\textbf{Case I}. Whenever $|b|\geq|s|+\frac{|a-s|}{12}$, 
we have the inequality, 
$$|ab-st|\geq|a-s|\cdot|b|-|b-t|\cdot|s|\geq\frac{|a-s|^2}{12}.$$
Finally, notice that $\|\pi_N(z)-\pi_N(\zeta)\|\geq|ab-st|$, so 
equation~(\ref{eq12}) holds in this particular case. 

\textbf{Case II}. If $|b|\leq|s|+\frac{|a-s|}{12}$, and there exists 
a natural $j$ such that $|a-\rho^j_Ns|$ is less than or equal to 
$\frac{|a-s|}{2}$, we also have,
$$|a-s|\leq|a-\rho^j_Ns|+2|s|\leq\frac{|a-s|}{2}+2|s|.$$ 
Consequently, $|s|\geq\frac{|a-s|}{4}$. On the other hand, we know that 
$\|z-\phi_j(\zeta)\|_\infty$ is equal to the maximum of $|a-\rho^j_Ns|$ 
and $|\rho^j_Nb-t|$. Recalling the hypotheses of Lemma~\ref{general} 
and this case (II), we have that $|a-\rho^j_Ns|<|a-s|$, and that 
$\|z-\phi_j(\zeta)\|_\infty$ is greater than or equal to 
$\|z-\zeta\|_\infty=|a-s|$. Hence, both $|\rho^j_Nb-t|\geq|a-s|$ and:
\begin{eqnarray*}
|ab-st|&\geq&|\rho^j_Nb-t|\cdot|s|-|a-\rho^j_Ns|\cdot|b|\\ &\geq&
\frac{|a-s|\cdot|s|}{2}-\frac{|a-s|^2}{24}\,\geq\,\frac{|a-s|^2}{12}.
\end{eqnarray*}
Notice that $|a-\rho^j_Ns|\cdot|b|$ is less than or equal to 
$\frac{|a-s|\cdot|s|}{2}+\frac{|a-s|^2}{24}$ because of the hypotheses of 
this case (II). We may conclude that equation~(\ref{eq12}) holds in this 
particular case as well, after recalling that $\|\pi_N(z)-\pi_N(\zeta)\|$ 
is greater than or equal to $|ab-st|$.

\textbf{Case III}. If $|b|\leq|s|+\frac{|a-s|}{12}$, and 
$|a-\rho^k_Ns|\geq\frac{|a-s|}{2}$ for every natural $k$, 
we automatically have the following inequality,
$$|a^N-s^N|=\prod^N_{k=1}|a-\rho^k_Ns|\geq\frac{|a-s|^N}{2^N}.$$

Finally, we know that $\|z-\zeta\|_{\infty}=|a-s|$, and that 
$\|\pi_N(z)-\pi_N(\zeta)\|$ is greater than or equal to $|a^N-s^N|$. 
The previous inequalities show that equation~(\ref{eq12}) holds for 
$\delta=N$. On the other hand, when $\delta=Ev(N)/2$, it is easy to 
deduce the existence of a subset $J$ of $\{1,2,\ldots,N\}$ composed 
of at least $N-\delta$ elements and which satisfies,
\begin{equation}\label{eq13}
|a-\rho_N^js|\geq\max\left\{|a|,|s|,\frac{|a-s|}{\sqrt{2}}\right\}
\quad\hbox{for each}\quad{j}\in{J}.
\end{equation}

The set $J$ can be built as follows. We may suppose, without lost of 
generality, that $a$ is real and $a\geq0$, for we only need to multiply 
both $a$ and $s$ by an appropriate complex number $\theta$ with 
$|\theta|=1$. Thus, the set $J$ is composed of all exponents 
$1\leq{j}\leq{N}$ which satisfy $\Re(\rho_N^js)\leq0$. It is easy to 
see that $|a-\rho_N^js|^2$ is greater than or equal to $|a|^2+|s|^2$, 
for every $j\in{J}$. Moreover, $J$ is composed of at least $N/2$ 
elements when $N$ is even, and of at least $\frac{N-1}{2}$ elements 
when $N$ is odd. Equation~(\ref{eq13}) follows automatically because 
$|a|^2$, $|s|^2$ and $\frac{|a-s|^2}{2}$ are all less than or equal 
to $|a|^2+|s|^2$. The hypotheses of this case (III), and 
equation~(\ref{eq13}), directly imply that~:
$$2|a-\rho_N^js|\geq|s|+\frac{|a-s|}{\sqrt{2}}\geq|b|,
\quad\forall{j}\in{J}.$$
Moreover, since $\frac{5}{3}>\frac{13}{12}\sqrt{2}$, and we are supposing 
from the beginning of this proof that $|a-s|\geq|b-t|$, we may also deduce 
the following inequality,
$$\frac{8|a-\rho_N^js|}{3}>|s|+\frac{|a-s|}{12}+|b-t|\geq|t|,
\quad\forall{j}\in{J}.$$
Finally, considering all the results presented in previous paragraphs, 
equation~(\ref{eq13}) and the hypotheses of this case (III), we can 
deduce the desired result,
$$|a^N-s^N|=\prod^N_{k=1}|a-\rho^k_Ns|\geq\frac{\|z,\zeta
\|^{N-\delta}_\infty|a-s|^\delta}{(8/3)^{N-\delta}\,2^\delta}$$
where $\delta=Ev(N)/2$, the norm $\|z-\zeta\|_\infty=|a-s|$ and 
$\|z,\zeta\|_\infty$ is the maximum of $|a|$, $|b|$, $|s|$ and $|t|$. 
We can conclude that equation~(\ref{eq12}) holds when $\delta$ is equal 
to $N$ and $Ev(N)/2$. 
\end{proof}

We are now in position to complete the proof of Theorem~\ref{main}. Notice 
that Lemma~\ref{general} automatically implies the following inequalities, 
whenever $z$ and $\zeta$ lie inside the compact ball $\overline{B_R}$, 
and $\delta=N$,
\begin{eqnarray*}
\|\pi_N(z)-\pi_N(\zeta)\|&\geq&\frac{\|z-\zeta\|^N_{\infty}}{2^N}\,
\min\left\{\frac{1}{3\,R^{N-2}}\,,1\right\}\\ 
&\geq&\frac{\|z-\zeta\|^N}{\sqrt{8}{^N}}\,
\min\left\{\frac{1}{3\,R^{N-2}}\,,1\right\}.
\end{eqnarray*}

\begin{proof}\textbf{(Theorem~\ref{main}, case $N\geq3$).} We shall follow 
step by step the proof of Theorem~\ref{main} (case $N=2$), presented in 
section two; so we shall only indicate the main differences. Let $g$ be a 
continuous solution to the equation $\dbar{g}=\pi^*\lambda$ on $B_R$ which 
satisfies the H\"older estimates given in equations~(\ref{eq2}) and 
(\ref{eq3}). Recalling the analysis done at the beginning of this section, 
we know there exists a continuous function $h$ defined on $\pi(B_R)$ such 
that $h\circ\pi$ is equal to $\frac{1}{N}\sum\phi_k^*g$. In particular, 
$\dbar{h}=\lambda$ on $\pi(B_R)$, and $\|h\|_{\pi_N(B_R)}$ is less than 
or equal to $\|g\|_{B_R}$. Moreover, working like in equations~(\ref{eq5}) 
and~(\ref{eq6}), we may deduce the existence of a finite positive constant 
$C_5(R)$ such that $\|\pi_N^*\lambda\|_{B_R}$ is less than or equal to 
$C_5(R)\|\lambda\|_{\pi_N(B_R)}$.

Given two points $x,w\in\pi(B_R)$, choose $z,\zeta\in{B_R}$ 
such that $x=\pi_N(z)$ and $w=\pi_N(\zeta)$. Since 
$\pi_N(\zeta)=\pi_N(\phi_k(\zeta))$ for every automorphism $\phi_k$ 
defined in~(\ref{eq11}), we can even choose $\zeta\in{B_R}$ so that 
$\|z-\zeta\|_{\infty}$ is less than or equal to 
$\|z-\phi_k(\zeta)\|_{\infty}$ for every $\phi_k$. A direct application 
of Lemma~\ref{general}, with $\delta=N$, yields the existence of a finite 
positive constant $C_6(R)$ such that $\|x-w\|$ is greater than or equal to 
$C_6(R)\|z-\zeta\|^{1/\beta}$. Recall that $\beta=1/Ev(N)$ and $N\geq3$. 
Thus, if $z$ and $\zeta$ are both inside the ball $B_{R/2}$, we may 
apply equation~(\ref{eq3}), in order to deduce that 
$\frac{|h(x)-h(w)|}{\|x-w\|^{\beta}}$ is less than or equal 
to $\frac{C_3(R)}{C_6^{\beta}(R)}\|\pi^*\lambda\|_{B_R}$.

On the other hand, suppose, without lost of generality, that $z$ is not 
inside the ball $B_{R/2}$. a direct application of Lemma~\ref{general}, 
with $\delta=\frac{Ev(N)}{2}=\frac{1}{2\beta}$, yields the existence of 
a finite positive constant $C_7(R)$ such that $\|x-w\|$ is greater than 
or equal to $C_7(R)\|z-\zeta\|^{\delta}$. Whence, equation~(\ref{eq2}) 
automatically implies that $\frac{|h(x)-h(w)|}{\|x-w\|^{\beta}}$ is less 
than or equal to $\frac{C_2(R)}{C_7^{\beta}(R)}\|\pi^*\lambda\|_{B_R}$ 
as well.

The analysis done in the previous paragraphs automatically implies 
the existence of a finite positive constant $C_1(R)$ such that 
equation~(\ref{eq1}) holds for every $N\geq3$.
\end{proof}

Finally, as we have already said at the end of section two, the proof 
of Theorem~\ref{main} works perfectly if we apply Theorem~\ref{Henkin} 
of Henkin and Leiterer, or any other integration kernel which produces 
estimates similar to those posed in equations~(\ref{eq2}) and (\ref{eq3}). 
For example, the hypotheses on $\lambda$ can be relaxed in 
Theorem~\ref{main}, to consider $(0,1)$-differential forms $\lambda$ which 
are bounded and continuous on $\pi(\overline{B_R})\setminus{K}$, for some 
compact set $K\subset{B_R}$ of zero-measure. Besides, the results 
presented in Theorem~\ref{main} hold as well, if we consider an arbitrary 
strictly pseudoconvex domain $D$, with smooth boundary and the origin in 
its interior, instead of the open ball $B_R$. In this case, the ball 
$B_{B/2}$ used in equation~(\ref{eq3}) of Theorem~\ref{Henkin} would be a 
sufficiently small ball $B_r$ whose closure is contained in the interior 
of $D$.

On the other hand, the work presented in this paper is strongly based 
on the existence of a branched finite covering $\pi_N$ from $\cc^2$ 
onto $X_N$, such that the inverse image of the singular point 
$\pi^{-1}_N(0)=\{0\}$ is a singleton. This property allows us to get the 
estimates presented in Lemmas~\ref{estimate} and~\ref{general}, which are 
essential for this paper. It is obvious to consider a \textit{blow-up} 
mapping $\eta:W\to{X_N}$ instead of the finite covering $\pi_N$. In any 
case, a blow-up is a $1$-covering everywhere, except at the singular 
point $0$. However, since the inverse image $\eta^{-1}(0)$ is not a 
singleton, and it is not even finite in general, we have strong problems 
for calculating a H\"older solution to the equation $\dbar{h}=\lambda$, 
unless we introduce stronger hypotheses. We finish this section by 
analysing the case of a blow-up.

\vspace{9pt}

\textbf{Remark.} Let $\Sigma$ be a variety with an isolated singularity 
at $\sigma_0\in\Sigma$, and $\eta:W\to\Sigma$ be a holomorphic 
\textit{blow-up} of $\Sigma$ at $\sigma_0$, such that $W$ is a smooth 
manifold. Given a $\dbar$-closed $(0,1)$-form $\lambda$ defined on 
$\Sigma$ minus $\sigma_0$, we automatically have that $\eta^*\lambda$ 
is also $\dbar$-closed on $W$ minus $\eta^{-1}(\sigma_0)$. Thus, 
suppose there exists a continuous solution $g:W\to\cc$ to the equation 
$\dbar{g}=\eta^*\lambda$. Since $\eta$ is a blow-up, we automatically 
have that $\eta^{-1}$ is well defined on $\Sigma\setminus\{\sigma_0\}$, 
and so $\lambda$ is equal to $\dbar(g\circ\eta^{-1})$ there.

Define $h:=g\circ\eta^{-1}$. Unless $g$ is constant on the inverse fibre 
$\eta^{-1}(\sigma_0)$, the function $h$ does not have a continuous 
extension to $\sigma_0$, and does not satisfy any H\"older condition in a 
neighbourhood of $\sigma_0$. Suppose there exists a pair of points $a$ and 
$b$ in $\eta^{-1}(\sigma_0)$ such that $g(a)\neq{g}(b)$. Besides, take 
$\{a_m\}$ and $\{b_m\}$ a pair of infinite sequences in 
$W\setminus\eta^{-1}(\sigma_0)$ which respectively converge to $a$ and 
$b$. Notice that both $\eta(a_m)$ and $\eta(b_m)$ converge to the same 
point $\sigma_0$. However, $g(a_m)$ and $g(b_m)$ converge to different 
points, for $g(a)\neq{g}(b)$. Hence, given a metric $\Delta$ on $\Sigma$, 
which defines the topology, we have that~:
$$\limsup_m\frac{|h\circ\eta(a_m)-h\circ\eta(b_m)|} 
{\Delta[\eta(a_m),\eta(b_m)]^\beta}=\infty,\quad\forall\beta>0.$$

That is, in order to introduce H\"older conditions on $h:=g\circ\eta^{-1}$, 
it is essential that the solution to equation $\dbar{g}=\eta^*\lambda$ is 
constant on the inverse fibre of the singular point $\eta^{-1}(\sigma_0)$.

\section{Surfaces with simple singularities of type $D_{N+2}$}

We finish this paper by solving the $\dbar$-equation on a neighbourhood 
of the origin in the subvariety $Y_N\subset\cc^3$, defined by the 
polynomial $y_1^2y_3+y_2^2=y_3^{N+1}$. The surface $Y_N$ has an isolated 
simple (rational double point) singularity of type $D_{N+2}$ at the origin, 
for any natural number $N\geq2$, \cite[p.~60]{Di}. We extend the results 
presented in theorem~\ref{main}, by introducing a branched $2$-covering 
defined from $X_{2N}:=[x_1x_2=x_3^{2N}]$ onto the surface $Y_N$. Consider 
the holomorphic mapping $\eta_2:X_{2N}\to{Y_N}$, and the pair of matrices 
$P$ and $Q$, given by the respective equations,
\begin{eqnarray}\label{eq21} 
&&\eta_2(x_1,x_2,x_3)=\left(\frac{x_1+x_2}{2}\,,
x_3\frac{x_1-x_2}{2i}\,,x_3^2\right),\\ 
\label{eq22}&&P=\left[\begin{array}{ccc} 0&1&0\\ 1&0&0\\ 0&0&-1
\end{array}\right],\quad{Q}=\frac{1}{2}\left[\begin{array}{ccc}
1&1&0\\ -i&i&0\\ 0&0&2\end{array}\right]. 
\end{eqnarray}

It is easy to see that $\eta_2(x)=\eta_2(Px)$ for every $x\in{X}_{2N}$. 
Moreover, $\eta_2$ is a branched $2$-covering, and the origin is the  
unique branch point, because the inverse image $\eta_2^{-1}(y)$ is a set 
of the form $\{x,Px\}$, for every $y\in{Y_N}$. For example, the inverse 
image of $(y_1,0,0)$ is composed of two points: $(2y_1,0,0)$ and 
$(0,2y_1,0)$. We have already defined a branched covering $\pi_{2N}$ from 
$\cc^2$ onto $X_{2N}$, so the composition $\eta_2\circ\pi_{2N}$ is indeed 
a covering from $\cc^2$ onto $Y_N$. The branch covering $\eta_2$ is a 
central part in the following result.

\begin{theorem}\label{DN}
Given an open ball $B_R\subset\cc^2$ of radius $R>0$ and centre in 
the origin, define the open set $E_R:=\eta_2(\pi_{2N}(B_R))$ in $Y_N$. 
There exists a finite positive constant $C_{11}(R)$ such that: for every 
continuous $(0,1)$-differential form $\aleph$ defined on the compact set 
$\overline{E_R}\subset{Y_N}$, and $\dbar$-closed on $E_R$, the equation 
$\dbar{f}=\aleph$ has a continuous solution $f$ on $E_R$ which also 
satisfies the following H\"older estimate, with $\beta_2=\frac{1}{4N}$,
\begin{equation}\label{eq23} 
\|f\|_{E_R}+\sup_{y,\xi\in{E_R}}\frac{|f(y)-f(\xi)|}
{\|y-\xi\|^{\beta_2}}\,\leq{C}_{11}(R)\|\aleph\|_{E_R}. 
\end{equation} 
\end{theorem}

The proof of this theorem follows exactly the same ideas and steps 
presented in the proof of Theorem~\ref{main} (case $N=2$), so we do 
not include it. Given a $(0,1)$-differential form $\aleph$ continuous 
on $\overline{E_R}\subset{Y_N}$, and $\dbar$-closed on $E_R$. We have 
that $\eta_2^*\aleph$ is also continuous on $\pi_{2N}(\overline{B_R})$, 
and $\dbar$-closed on $\pi_{2N}(B_R)$. Therefore, we can apply 
Theorem~\ref{main}, in order to obtain a continuous solution $h$ to the 
differential equation $\dbar{h}=\eta_2^*\aleph$, which also satisfies 
the H\"older conditions given in equation~(\ref{eq1}). There exists a 
continuous function $f$ on $E_R$ such that $\eta_2^*f$ is equal to 
$\frac{h+\psi^*h}{2}$, and so $\dbar{f}=\aleph$, as we wanted. Finally, 
inequality~(\ref{eq23}) follows from equation~(\ref{eq1}), after noticing 
that there exists a pair of finite positive constants $C_8(R)$ and 
$C_9(R)$ such that~:
$$\|f\|_{E_R}\leq\|h\|_{\pi_{2N}(B_R)},\quad
\|\eta_2^*\aleph\|_{\pi_{2N}(B_R)}\leq{C}_8(R)\|\aleph\|_{E_R},$$ 
and $\frac{|f(y)-f(\xi)|}{\|y-\xi\|^{\beta_2}}$ is also less than or equal 
to $C_9(R)\|\eta_2^*\aleph\|_{\pi_{2N}(B_R)}$ for every $y$ and $\xi$ in 
$E_R$. We obviously need an estimate of $\|x-w\|^2$, with respect to the 
projections $\|\eta_2(x)-\eta_2(w)\|$, in order to show that the previous 
inequality above holds. This estimate is presented in the following
Lemma~\ref{final}. In conclusion, the proof of Theorem~\ref{DN} follows 
the same ideas and steps of the proof of Theorem~\ref{main} (case $N=2$), 
we only need to apply Theorem~\ref{main} instead of Theorem~\ref{Henkin}, 
and the following Lemma~\ref{final} instead of Lemma~\ref{estimate}.

\begin{lemma}\label{final}
Let $x$ and $w$ be two points in $X_{2N}$ whose norms $\|x\|$ and $\|w\|$ 
are both less than or equal to a finite constant $\rho>0$. If the distance 
$\|Q(w-Px)\|$ is greater than or equal to $\|Q(w-x)\|$, for the matrices 
$P$ and $Q$ defined in~(\ref{eq22}), then, the following inequality holds. 
\begin{eqnarray} 
\label{eq24}\|\eta_2(x)-\eta_2(w)\|&\geq&C_{12}(\rho)\|x-w\|^2,\\ 
\nonumber\hbox{where}\quad{C}_{12}(\rho)&=&\frac{1}{80}\min 
\left\{4,\frac{5}{3\rho}\,,\frac{1}{\rho^{2N-2}N}\right\}. 
\end{eqnarray} 
\end{lemma}

\begin{proof} Introducing the new variables $(a,b,c):=Qx$ and 
$(s,t,u):=Qw$, we have that $a^2+b^2=c^{2N}$ and $QPx=(a,-b,-c)$ 
for every $x\in{X}_{2N}$. Moreover,
\begin{equation}\label{eq25}
\|\eta_2(x)-\eta_2(w)\|^2=|a-s|^2+|bc-tu|^2+|c^2-u^2|^2.
\end{equation}

A main step in this proof is to shown that the following inequality holds,
\begin{equation}\label{eq26}
\|\eta_2(x)-\eta_2(w)\|\geq16C_{12}(\rho)\|\pi_2(b,c)-\pi_2(t,u)\|,
\end{equation}
where $\pi_2(b,c)=(b^2,c^2,bc)$ was defined in the introduction of this 
paper, and $C_{12}(\rho)$ is given in equation~(\ref{eq24}) above. We know 
that $\|Q(w-Px)\|$ is greater than or equal to $\|Q(w-x)\|$, according to 
the hypotheses of this lemma, so it is easy to deduce that $\|t+b,u+c\|$ 
is also greater than or equal to $\|t-b,u-c\|$, because $QPx$ is equal 
to $(a,-b,-c)$. Therefore, if equation~(\ref{eq26}) holds, a direct 
application of Lemma~\ref{estimate} yields,
\begin{equation}\label{eq27}
\|\eta_2(x)-\eta_2(w)\|\geq8C_{12}(\rho)(|b-t|^2+|c-u|^2).
\end{equation}

On the other hand, we can easily calculate the following upper bound
for $|a|$,
\begin{equation}\label{eq28}
|a|\leq\|Qx\|\leq\|x\|\leq\rho.
\end{equation}
A similar upper bound $|s|\leq\rho$ holds as well. Hence, recalling 
equation~(\ref{eq25}), we have that $\|\eta_2(x)-\eta_2(w)\|$ is 
greater than or equal to $|a-s|\geq\frac{|a-s|^2}{2\rho}$. Adding together 
the inequality presented in previous statement and equation~(\ref{eq27}) 
yields the desired result, notice that $\frac{1}{2\rho}>8C_{12}(\rho)$ 
and $2\|\xi\|\geq\|Q^{-1}\xi\|$ for $\xi\in\cc^3$,
\begin{eqnarray*}
2\|\eta_2(x)-\eta_2(w)\|&\geq&8C_{12}(\rho)\|Q(x-w)\|^2\\
&\geq&2C_{12}(\rho)\|w-x\|^2.
\end{eqnarray*}

We may then conclude that inequality~(\ref{eq24}) holds, as we wanted. 
We only need to prove that equation~(\ref{eq26}) is always satisfied, 
in order to finish our calculations; and we are doing to prove that by 
considering two complementary cases.

\textbf{Case I} Whenever $3|c^2-u^2|$ is greater than or equal 
to $\frac{|b^2-t^2|}{\rho^{2N-2}N}$, the following inequality holds,
$$|c^2-u^2|^2\geq\frac{16\,|c^2-u^2|^2}{25}+
\frac{|b^2-t^2|^2}{(5\rho^{2N-2}N)^2}.$$
Thus, in this particular case, inequality~(\ref{eq26}) follows 
directly from equation~(\ref{eq25}), because $\frac{1}{5\rho^{2N-2}N}$ 
and $4/5$ are both greater than or equal to $16C_{12}(\rho)$.

\textbf{Case II} Whenever $\frac{|b^2-t^2|}{\rho^{2N-2}N}$ is greater 
than or equal to $3|c^2-u^2|$, we proceed as follows. The absolute 
values $|a|$ and $|c|$ are both bounded by $\|Qx\|\leq\rho$, according 
to equation~(\ref{eq28}); the same upper bound can be calculated for 
$|s|$ and $|u|$. Whence, the following series of inequalities hold~:
\begin{eqnarray*}
2|b^2-t^2|/3&\leq&|b^2-t^2|-\rho^{2N-2}N|c^2-u^2|\\
&\leq&|b^2-t^2|-|c^{2N}-u^{2N}|\\
&\leq&|a^2-s^2|\;\leq\;2\rho|a-s|.
\end{eqnarray*}
Recall that $a^2+b^2=c^{2N}$ and that $(\xi^N-1)$ is equal to 
the product of $(\xi-1)$ times the sum $\sum_{k=0}^{N-1}\xi^k$. 
Inequality~(\ref{eq26}) follows then from equation~(\ref{eq25}), 
after noticing that $16C_{12}(\rho)|b^2-t^2|$ is less than or equal 
to $|a-s|$, and obviously, $16C_{12}(\rho)$ is also less than one.
\end{proof}

\bibliographystyle{plain}

\end{document}